\documentclass[english,conference]{IEEEtran}
\usepackage[T1]{fontenc}
\usepackage[latin9]{inputenc}
\usepackage{geometry}
\usepackage{amsthm}
\usepackage{amsmath}
\usepackage{amssymb}

\makeatletter
\theoremstyle{plain}
\newtheorem{thm}{\protect\theoremname}
\theoremstyle{definition}
\newtheorem{defn}[thm]{\protect\definitionname}
\theoremstyle{plain}
\newtheorem{lem}[thm]{\protect\lemmaname}
\theoremstyle{plain}
\newtheorem{prop}[thm]{\protect\propositionname}
\theoremstyle{definition}
\newtheorem{example}[thm]{\protect\examplename}
\theoremstyle{remark}
\newtheorem{rem}[thm]{\protect\remarkname}

\usepackage{cite}
\author{\IEEEauthorblockN{Farzad Farnoud (Hassanzadeh), Behrouz Touri, and Olgica Milenkovic}
\authorblockA{University of Illinois, Urbana-Champaign, IL}
E-mail: \{hassanz1,touri1,milenkov\}@illinois.edu}

\makeatother

\usepackage{babel}
\providecommand{\definitionname}{Definition}
\providecommand{\examplename}{Example}
\providecommand{\lemmaname}{Lemma}
\providecommand{\propositionname}{Proposition}
\providecommand{\remarkname}{Remark}
\providecommand{\theoremname}{Theorem}

\begin{document}
\global\long\def\dist{\mathsf{d}}
\global\long\def\kdist{\mathsf{d}_{\tau}}
\global\long\def\Kdist{\mathsf{d}_{K}}
\global\long\def\Tdist{\mathsf{d}_{T}}
\global\long\def\define{:\,=}
\global\long\def\dash{\mbox{--}}
\global\long\def\adjset{A}
\global\long\def\lineset{L}

\title{\vspace{10pt}Novel Distance Measures for Vote Aggregation}
\maketitle
\begin{abstract}
We consider the problem of rank aggregation based on new distance
measures derived through axiomatic approaches and based on score-based
methods. In the first scenario, we derive novel distance measures
that allow for discriminating between the ranking process of highest
and lowest ranked elements in the list. These distance functions represent
weighted versions of Kendall's $\tau$ measure and may be computed
efficiently in polynomial time. Furthermore, we describe how such
axiomatic approaches may be extended to the study of score-based aggregation
and present the first analysis of distributed vote aggregation over
networks.
\end{abstract}

\section{Introduction}

Rank aggregation is a classical problem frequently encountered in
social sciences, web search and Internet service analysis, expert
opinion and voting theory \cite{kemeney1959mathematics,cook1985ordinal,dwork2001rank,dwork2001rank-web,sculley2007rank,schalekamp2009rank,kumar2010gdr}.
The problem can be succinctly described as follows: a set of ``voters''
or ``experts'' is presented with a set of distinguishable entities
(objects, individuals, movies), typically represented by the set $\{1,2,\cdots,n\}$.
The voters' task is to arrange the entities in decreasing order of
preference and pass on their ordered lists to an aggregator. The aggregator
outputs a single preference list used as a representative of all voters.
Hence, one has to be able to adequately measure the quality of representation
made by a vote aggregator. Two distinct analytical rank aggregation
methods were proposed so far, namely, \emph{distance-based} methods
and \emph{score-(position-)based} methods. In the first case, the
quality of the aggregate is measured via a distance function that
describes how close the aggregate is to each individual vote. In the
second case, the aggregate is obtained by computing a score for each
ranked entity and then arranging the entities based on their score.
Well known distance measures include Kendall\textquoteright{}s $\tau$
and Spearman\textquoteright{}s Footrule \cite{diaconis1988group}.

The goal of this work is to propose two novel research directions
in rank aggregation: one, which builds upon the existing work of distance-based
aggregation, but expands the scope and applicability of vote-distances;
and another, which sets the stage for analyzing score-based vote aggregations
over networks. The results presented in the paper include a new set
of voting-fairness axioms that lead to distance measures previously
unknown in literature, as well as an analysis of consensus in distributed
score-based voting systems.

Our work on aggregation distance analysis is motivated by the following
observations: a) in many applications, the top of the ranking is more
important than the bottom and so changes to the top of the list must
result in a more significant change in the aggregate ranking than
changes to the bottom of the list; b) ranked entities may have different
degrees of similarity and often the goal is to find the most diverse,
yet highest ranked entities. Hence, swapping elements that are similar
should be penalized less than swapping those that are not. To the
best of the authors' knowledge, the work of Sculley \cite{sculley2007rank}
represents the only method proposed so far for handling similarity
in rank aggregation. Sculley presents an aggregation method, based
on the use of Markov chains first introduced by Dwork et. al., with
the goal of assigning similar ranks to similar items. A handful of
results are known for rank aggregation distances that address the
problem of positional relevance, i.e. the significance of the top
versus the bottom of the ranking \cite{kumar2010gdr}. In this context,
we introduce the notions of weighted Kendall distance and weighted
Cayley distance, both capable of addressing the top versus bottom
ranking issue, and provide axiomatic characterizations for these distance
measures.

The work on vote aggregation over networks considers the issue of
reaching consensus about the aggregate ranking in an arbitrary network,
either through local interactions or based on a gossip algorithms.
The assumption is that voters are connected through a social network
that allows them to adjust their votes based on the opinions of their
neighbors or randomly chosen network nodes, or even based on exogenous
opinions. For a special type of score-based scheme -- Borda's rule
-- we show that convergence to a vote consensus occurs and we determine
the rate of convergence. The analysis of rank aggregation over networks
for distance-based aggregation rules, and in particular for Kendall's
$\tau$ and the weighted Kendall distance, is postponed to the full
version of the paper.

The paper is organized as follows. An overview of relevant concepts,
definitions, and terminology is presented in Section II. Weighted
Kendall distance measures and extensions thereof, as well their axiomatic
definitions, are presented in Sections III and IV. Section V is devoted
to the analysis of gossip algorithms for rank aggregation.

\section{Preliminaries}

Suppose one is given a set $\Sigma=\{\sigma_{1},\sigma_{2},\cdots,\sigma_{m}\}$
of rankings, where each ranking $\sigma$ represents a permutation
in $\mathbb{S}_{n}$, the symmetric group of order $n$.

Given a distance function $\dist$ over the permutations in $\mathbb{S}_{n}$,
the distance-based aggregation problem can be stated as 
\[
\min_{\pi\in\mathbb{S}_{n}}\sum_{i=1}^{m}\dist(\pi,\sigma_{i}).
\]
 In words, the goal is to find a ranking $\pi$ with minimum cumulative
distance from $\Sigma$. Clearly, the choice of the distance function
$\dist$ is an important feature for all distance-based rank aggregation
methods. Many distance measures in use were derived by starting from
a reasonable set of axioms and then showing that the given distance
measure is a unique solution under the given set of axioms%
\footnote{This is to be contrasted with the celebrated Arrow's impossibility
theorem \cite{arrow1963social}.%
}. A distance function derived in this manner is Kendall's $\tau$
distance, based on Kemeny's axioms \cite{kemeney1959mathematics}.

On the other hand, score-based methods are centered around aggregators
that assign scores to objects based on their positions in the rankings
of $\Sigma$. Objects are then sorted according to their scores to
obtain the aggregate ranking. One of the best known rules in this
family is Borda's aggregation rule, introduced by Jean-Charles de
Borda \cite{borda1784} wherein, for each ranking $\sigma_{i}$, object
$j$ receives score $b_{i}^{j}=\sigma_{i}^{-1}(j)$. The average score
of object $j$ is $\bar{b}^{j}=\frac{1}{m}\sum_{i=1}^{m}b_{i}^{j}$.
The aggregate ranking is obtained by assigning the highest rank to
the object with the lowest average score, the second highest rank
to the object with the second lowest average score and so on. Borda's
method also has an axiomatic underpinning: in the context of social
choice functions, Young \cite{young1974-an-axiomatization} presented
a set of axioms that showed that Borda's rule is the unique social
choice function that satisfies the given axioms. A social choice function
is a rule indicating a set of winners when votes are given as rankings.
Note that although similar, a social choice function differs from
an aggregation rule; while a social choice function returns a set
of winners, an aggregation rule ranks all objects.

In what follows, we introduce the notation used throughout the paper
and provide a novel proof for the uniqueness of Kendall's $\tau$
distance function for a set of reduced Kemeny axioms \cite{kemeney1959mathematics}. 

Let $e=12\cdots n$ denote the identity permutation (ranking).
\begin{defn}
A transposition of two elements $a,b\in[n]$ in a permutation $\pi$
is the swap of elements in positions $a$ and $b$, and is denoted
by $(a\, b)$. In general, we reserve the notation $\tau$ for an
arbitrary transposition and when there is no confusion, we consider
a transposition to be a permutation itself. If $|a-b|=1$, the transposition
is referred to as an \emph{adjacent transposition.}

It is well known that any permutation may be reduced to $e$ via transpositions
or adjacent transpositions. The former process is referred to as sorting,
while the later is known as sorting with adjacent transpositions.
The smallest number of adjacent transpositions needed to sort a permutation
$\pi$ is known as the inversion number of the permutation. Equivalently,
the corresponding distance $\dist(e,\pi)$ is known as the Kendall's
$\tau$ distance. The Kendall's $\tau$ can be computed in time $O(n^{2}).$

We also find the following set useful in our analysis, 
\begin{multline*}
\adjset(\pi,\sigma)=\bigl\{\left(\tau_{1},\cdots,\tau_{m}\right):\\
m\in\mathbb{N},\sigma=\pi\tau_{1}\cdots\tau_{m},\tau_{i}=\left(a_{i}\ a_{i}+1\right),i\in[m]\bigr\}
\end{multline*}
i.e., the set of all sequences of adjacent transpositions that transform
$\pi$ into $\sigma.$

For a ranking $\pi\in\mathbb{S}_{n}$ and $a,b\in[n]$, $\pi$ is
said to \emph{rank} $a$ \emph{before} $b$ if $\pi^{-1}(a)<\pi^{-1}(b)$.
We denote this relationship as $a<_{\pi}b$. Two rankings $\pi$ and
$\sigma$ \emph{agree} on a pair $\{a,b\}$ of elements if both rank
$a$ before $b$ or both rank $b$ before $a$. Furthermore, the two
rankings $\pi$ and $\sigma$ \emph{disagree} on the pair $\{a,b\}$
if one ranks $a$ before $b$ and the other ranks $b$ before $a$.

For example, consider $\pi=1234$ and $\sigma=4213$. We have that
$4<_{\sigma}1$ and that $\pi$ and $\sigma$ agree for $\{2,3\}$
but disagree for $\{1,2\}$.
\end{defn}

\begin{defn}
\label{def:between} A ranking $\omega$ is said to be \emph{between
}two rankings $\pi$ and $\sigma$, denoted by $\pi\dash\omega\dash\sigma$,
if for each pair of elements $\{a,b\}$, $\omega$ either agrees with
$\pi$ or $\sigma$ (or both). The rankings $\pi_{0},\cdots,\pi_{m}$
are said to be on a line, denoted by $\pi_{0}\dash\pi_{1}\dash\cdots\dash\pi_{m}$,
if for every $i,j,$ and $k$ for which $0\le i<j<k\le m$, we have
$\pi_{i}\dash\pi_{j}\dash\pi_{k}$.
\end{defn}
The basis of our subsequent analysis is the following set of axioms
required for a rank aggregation measure, first introduced by Kemeny
\cite{kemeney1959mathematics}:

\textbf{Axioms I}
\begin{enumerate}
\item $\dist$ is a metric.
\item $\dist$ is left-invariant, i.e. $\dist(\sigma\pi,\sigma\omega)=\dist(\pi,\omega)$,
for any $\pi,\sigma,\omega\in\mathbb{S}_{n}$. In words, relabeling
of objects should not change the distance between permutations.
\item For any $\pi,\sigma,$ and $\omega$, $\dist(\pi,\sigma)=\dist(\pi,\omega)+\dist(\omega,\sigma)$
if and only if $\omega$ is between $\pi$ and $\sigma$. This axiom
may be viewed through a geometric lens: the triangle inequality has
to be satisfied for all points that lie on a ``straight line'' between
$\pi$ and $\sigma$.
\item The smallest positive distance is one. This axiom is only used for
normalization purposes. 
\end{enumerate}
Kemeny's original exposition included a fifth axiom which we restate
for completeness: If two rankings $\pi$ and $\sigma$ agree except
for a segment of $k$ elements, the position of the segment within
the ranking is not important. Here, a segment represents a set of
objects that are ranked consecutively - i.e., a substring of the permutation.
As an example, this axiom implies that 
\[
\dist(123\underbrace{456},123\underbrace{654})=\dist(1\underbrace{456}23,1\underbrace{654}23)
\]
 where the segment is denoted by braces. This axiom is redundant since
an equally strong statement follows from the other four axioms, as
we demonstrate below. Our alternative proof of Kemeny's result also
reveals a simple method for generalizing the axioms in order to arrive
at weighted distance measures.
\begin{lem}
\label{lem:on-the-line}For any $\dist$ that satisfies Axioms I,
and for any set of permutations $\pi_{0},\cdots,\pi_{m}$ such that
$\pi_{0}\dash\pi_{1}\dash\cdots\dash\pi_{m}$, one has 
\[
\dist(\pi_{0},\pi_{m})=\sum_{k=1}^{m}\dist(\pi_{k-1},\pi_{k}).
\]
\end{lem}
\begin{IEEEproof}
The lemma follows from Axiom I.3 by induction.\end{IEEEproof}
\begin{lem}
\label{lem:two-ways}For any $\dist$ that satisfies Axioms I, we
have that $\dist\left(\left(i\ i+1\right),e\right)=\dist\left((12),e\right),i\in[n-1]$.\end{lem}
\begin{IEEEproof}
We show that $\dist\left(\left(23\right),e\right)=\dist\left(\left(12\right),e\right)$.
Repeating the same argument used for proving this special case gives
$\dist\left(\left(i\ i+1\right),e\right)=\dist\left(\left(i-1\ i\right),e\right)=\cdots=\dist\left(\left(12\right),e\right)$.

To show that $\dist\left(\left(23\right),e\right)=\dist\left(\left(12\right),e\right)$,
we evaluate $\dist(\pi,e)$ in two ways, where we choose $\pi=32145\cdots n.$

On the one hand, note that $\pi\dash\omega\dash\eta\dash e$ where
$\omega=\pi(12)=23145\cdots n$ and $\eta=\omega(23)=21345\cdots n$.
As a result, 
\begin{align}
\dist(\pi,e) & \stackrel{\mathsf{(a)}}{=}\dist(\pi,\omega)+\dist(\omega,\eta)+\dist(\eta,e)\nonumber \\
 & =\dist(\omega^{-1}\pi,e)+\dist(\eta^{-1}\omega,e)+\dist(\eta,e)\nonumber \\
 & =\dist((12),e)+\dist((23),e)+\dist((12),e)\label{eq:one-hand}
\end{align}
where $(\mathsf{a})$ follows from Lemma \ref{lem:on-the-line}.

On the other hand, note that $\pi\dash\alpha\dash\beta\dash e$ where
$\alpha=\pi(23)=31245\cdots n$ and $\beta=\alpha(12)=13245\cdots n$.
For this case,
\begin{align}
\dist(\pi,e) & =\dist(\pi,\alpha)+\dist(\alpha,\beta)+\dist(\beta,e)\nonumber \\
 & =\dist(\alpha^{-1}\pi,e)+\dist(\beta^{-1}\alpha,e)+\dist(\beta,e)\nonumber \\
 & =\dist((23),e)+\dist((12),e)+\dist((23),e).\label{eq:other-hand}
\end{align}

Expressions (\ref{eq:one-hand}) and (\ref{eq:other-hand}) imply
that $\dist\left(\left(23\right),e\right)=\dist\left(\left(12\right),e\right)$.\end{IEEEproof}
\begin{lem}
\label{lem:gamma-1}For any $\dist$ that satisfies Axioms I, $\dist(\gamma,e)$
equals the minimum number of adjacent transpositions required to transform
$\gamma$ into $e$.\end{lem}
\begin{IEEEproof}
Let 
\[
\lineset(\pi,\sigma)=\left\{ \left(\tau_{1},\cdots,\tau_{m}\right)\in\adjset(\pi,\sigma):\pi\dash\pi\tau_{1}\dash\pi\tau_{1}\tau_{2}\dash\cdots\dash\sigma\right\} 
\]
be the subset of $\adjset(\pi,\sigma)$ consisting of sequences of
transpositions that transform $\pi$ to $\sigma$ by passing through
a line. Let $m$ be the minimum number of adjacent transpositions
that transform $\gamma$ into $e$. Furthermore, let $(\tau_{1},\tau_{2},\cdots,\tau_{m})\in A(\gamma,e)$
and define $\gamma_{i}=\gamma\tau_{1}\cdots\tau_{i},i=0,\cdots,m,$
with $\gamma_{0}=\gamma$ and $\gamma_{m}=e$. 

First, we show that $\gamma_{0}\dash\gamma_{1}\dash\cdots\dash\gamma_{m}$,
that is, $(\tau_{1},\tau_{2},\cdots,\tau_{m})\in L(\gamma,e)$. Suppose
this were not the case. Then, there would exist $i<j<k$ such that
$\gamma_{i},\gamma_{j},$ and $\gamma_{k}$ are not on a line, and
thus, there would exists a pair $\{r,s\}$ for which $\gamma_{j}$
disagrees with both $\gamma_{i}$ and $\gamma_{k}$. Hence, there
would be two transpositions, $\tau_{i'}$ and $\tau_{j'}$, with $i<i'\le j$
and $j<j'\le k$ that swap $r$ and $s$. We could in this case remove
$\tau_{i'}$ and $\tau_{j'}$ from $\left(\tau_{0},\cdots,\tau_{m}\right)$
to obtain $\left(\tau_{0},\cdots,\tau_{i'-1},\tau_{i'+1},\cdots,\tau_{j'-1},\tau_{j'+1},\tau_{m}\right)\in A(\gamma,e)$
with length $m-2$. This contradicts the optimality of the choice
of $m$. Hence, $(\tau_{1},\tau_{2},\cdots,\tau_{m})\in L(\gamma,e)$.
Then Lemma \ref{lem:on-the-line} implies that 
\begin{equation}
\dist(\gamma,e)=\sum_{i=1}^{m}\dist(\tau_{i},e).\label{eq:A-1}
\end{equation}

From (\ref{eq:A-1}), it is clear that the minimum positive distance
from the identity is obtained by some adjacent transpositions. But,
Lemma \ref{lem:two-ways} states that all adjacent transpositions
have the same distance from the identity. Hence, from Axiom I.4, we
have $\dist(\tau,e)=1$ for all adjacent transposition $\tau$. This
observation completes the proof of the lemma, since it implies that
\[
\dist(\gamma,e)=\sum_{i=1}^{m}\dist(\tau_{i},e)=\sum_{i=1}^{m}1=m.
\]
\end{IEEEproof}
\begin{lem}
\label{lem:inversions}For any $\dist$ that satisfies Axioms I, we
have

\[
\dist(\pi,\sigma)=\min\left\{ m:(\tau_{1},\cdots,\tau_{m})\in A(\pi,\sigma)\right\} .
\]
That is, $\dist(\pi,\sigma)$ equals the minimum number of adjacent
transpositions required to transform $\pi$ into $\sigma$.\end{lem}
\begin{IEEEproof}
We have $(\tau_{1},\cdots,\tau_{m})\in A(\pi,\sigma)$ if and only
if $(\tau_{1},\cdots,\tau_{m})\in A(\sigma^{-1}\pi,e)$. Furthermore,
left-invariance of $\dist$ implies that $\dist(\pi,\sigma)=\dist(\sigma^{-1}\pi,e)$.
Hence,
\begin{align*}
\dist(\pi,\sigma) & =\dist(\sigma^{-1}\pi,e)\\
 & =\min\left\{ m:(\tau_{1},\cdots,\tau_{m})\in A(\sigma^{-1}\pi,e)\right\} \\
 & =\min\left\{ m:(\tau_{1},\cdots,\tau_{m})\in A(\pi,\sigma)\right\} 
\end{align*}
where the second equality follows from Lemma \ref{lem:gamma-1}.\end{IEEEproof}
\begin{thm}
The unique distance $\dist$ that satisfies Axioms I is 
\[
\kdist(\pi,\sigma)=\min\left\{ m:(\tau_{1},\cdots,\tau_{m})\in A(\pi,\sigma)\right\} .
\]
\end{thm}
\begin{IEEEproof}
The fact that $\kdist$ satisfies Axioms I can be easily verified.
Uniqueness follows from Lemma \ref{lem:inversions}.
\end{IEEEproof}

\section{Weighted Kendall Distance}

Our proof of the uniqueness of Kendall's $\tau$ distance under Axioms
I reveals an important insight: Kendall's measure arises due to the
fact that adjacent transpositions have uniform costs, which is a consequence
of the betweenness property of one of the axioms. If one had a ranking
problem in which costs of transpositions either depended on the elements
involved or their locations, the uniformity assumption had to be changed.
As we show below, a way to achieve this goal is to redefine the axioms
in terms of the betweenness property.

\textbf{Axioms II}
\begin{enumerate}
\item $\dist$ is a pseudo-metric, i.e. a generalized metric in which two
distinct points may be at zero distance.
\item $\dist$ is left-invariant.
\item For any $\pi,\sigma$ disagreeing for more than one pair of elements,
there exists a $\omega$ such that $\dist(\pi,\sigma)=\dist(\pi,\omega)+\dist(\omega,\sigma)$.\end{enumerate}
\begin{lem}
\label{lem:WKD}For any distance $\dist$ that satisfies Axioms II,
and for distinct $\pi$ and $\sigma$, we have
\[
\dist(\pi,\sigma)=\min_{\left(\tau_{0},\cdots,\tau_{m}\right)\in\adjset(\pi,\sigma)}\sum_{i=1}^{m}\dist(\tau_{i},e).
\]
\end{lem}
\begin{IEEEproof}
First, suppose that $\pi$ and $\sigma$ disagree on one pair of elements.
Then, we have $\sigma=\pi(a\ a+1)$ for some $a\in[n-1]$. For each
$\left(\tau_{0},\cdots,\tau_{m}\right)\in\adjset(\pi,\sigma)$, there
exists an index $j$ such that $\tau_{j}=(a\ a+1)$ and thus 
\[
\sum_{i=1}^{m}\dist(\tau_{i},e)\ge\dist(\tau_{j},e)=\dist(\left(a\ a+1\right),e)
\]
implying
\begin{equation}
\min_{\left(\tau_{0},\cdots,\tau_{m}\right)\in A(\pi,\sigma)}\sum_{i=1}^{m}\dist(\tau_{i},e)\ge\dist(\left(a\ a+1\right),e).\label{eq:lq}
\end{equation}
 On the other hand, since $\left((a\ a+1)\right)\in\adjset(\pi,\sigma)$,
\begin{equation}
\min_{\left(\tau_{0},\cdots,\tau_{m}\right)\in A(\pi,\sigma)}\sum_{i=1}^{m}\dist(\tau_{i},e)\le\dist(\left(a\ a+1\right),e).\label{eq:hq}
\end{equation}
From (\ref{eq:lq}) and (\ref{eq:hq}), 
\[
\min_{\left(\tau_{0},\cdots,\tau_{m}\right)\in A(\pi,\sigma)}\sum_{i=1}^{m}\dist(\tau_{i},e)=\dist(\left(a\ a+1\right),e)=\dist(\pi,\sigma)
\]
where the last equality follows from the left-invariance of $\dist$.

Next, suppose $\pi$ and $\sigma$ disagree for more than one pair
of elements. A sequential application of Axiom II.3 implies that 
\[
\dist(\pi,\sigma)=\min_{\left(\tau_{0},\cdots,\tau_{m}\right)\in A(\pi,\sigma)}\sum_{i=1}^{m}\dist(\tau_{i},e),
\]

which proves the claimed result.\end{IEEEproof}
\begin{defn}
A distance $\dist$ is termed a \emph{weighted Kendall distance }if
there is a nonnegative \emph{weight function} $\varphi$ over the
set of adjacent transpositions such that 
\[
\dist(\pi,\sigma)=\min_{\left(\tau_{0},\cdots,\tau_{m}\right)\in\adjset(\pi,\sigma)}\sum_{i=1}^{m}\varphi_{\tau_{i}}
\]
where $\varphi_{\tau}$ is the weight assigned to transposition $\tau$
by $\varphi$.
\end{defn}
Note that a weighted Kendall distance is completely determined by
its weight function $\varphi$.
\begin{thm}
\label{thm:WKD}A distance $\dist$ satisfies Axioms II if and only
if it is a weighted Kendall distance.\end{thm}
\begin{IEEEproof}
It follows immediately from Lemma \ref{lem:WKD} that a distance $\dist$
satisfying Axioms II is a weighted Kendall distance by letting 
\[
\varphi_{\tau}=\dist(\tau,e)
\]
for every adjacent transposition $\tau$.

The proof of the converse is omitted since it is easy to verify that
a weighted Kendall distance satisfies Axioms II.
\end{IEEEproof}

The weighted Kendall distance provides a natural solution for issues
related to the importance of the top-ranked candidates. Due to space
limitations, we refer the reader interested in other applications
of weighted distances to our recent work \cite{farnoud2012sorting}.

\subsection*{Computing the Weighted Kendall Distance}

\global\long\def\wtfn{\varphi}
\global\long\def\wt#1{\varphi_{#1}}

Computing the weighted Kendall distance between two rankings for general
weight functions is not a task as straightforward as computing the
Kendall's $\tau$ distance. However, in what follows, we show that
for an important class of weight functions -- termed ``monotonic''
weight functions -- the weighted Kendall distance can be computed
efficiently. 
\begin{defn}
A weight function $\wtfn:\mathbb{A}_{n}\to\mathbb{R}^{+}$, where
$\mathbb{A}_{n}$ is the set of adjacent transpositions in $\mathbb{S}_{n}$,
is decreasing if $i>j$ implies that $\wt{(i\ i+1)}\le\wt{(j\ j+1)}.$
Increasing weight functions are defined similarly.
\end{defn}
Decreasing weight functions are important as they can be used to model
the significance of the top of the ranking by assigning higher weights
to transpositions at the top of the list.

Suppose a transformation $\tau=\left(\tau_{1},\cdots,\tau_{m}\right)$
of length $m$ transforms $\pi$ into $\sigma$. The transformation
may be viewed as a sequence of moves of elements indexed by $i$,
$i=1,\ldots,m,$ from position $\pi^{-1}(i)$ to position $\sigma^{-1}(i)$.
Let the \emph{walk} along which element $i$ is moved by transformation
$\tau$ be denoted by $p^{i,\tau}=\left(p_{1}^{i,\tau},\cdots,p_{\left|p^{i,\tau}\right|+1}^{i,\tau}\right)$
where $\left|p^{i,\tau}\right|$ is the length of the walk $p^{i,\tau}$. 

We investigate the lengths of the walks $p^{i,\tau},i\in[n].$ Let
$\mathcal{I}_{i}(\pi,\sigma)$ be the set consisting of elements $j\in[n]$
such that $\pi$ and $\sigma$ disagree on the pair $\{i,j\}$. Furthermore,
let $I_{i}(\pi,\sigma)=\left|\mathcal{I}_{i}(\pi,\sigma)\right|$.
In the transformation $\tau$, all elements of $\mathcal{I}_{i}(\pi,\sigma)$
must be swapped with $i$ by some $\tau_{k},k\in[m]$. Each such swap
contributes length one to the walk $p^{i,\tau}$ and thus, $\left|p^{i,\tau}\right|\ge I_{i}(\pi,\sigma)$. 

It is easy to see that 
\[
\dist_{\wtfn}(\pi,\sigma)=\min_{\tau\in A(\pi,\sigma)}\sum_{i=1}^{n}\frac{1}{2}\sum_{j=1}^{\left|p^{i,\tau}\right|}\wt{\left(p_{j}^{i,\tau}p_{j+1}^{i,\tau}\right)}.
\]
Considering individual walks, we may thus write
\begin{equation}
\dist_{\wtfn}(\pi,\sigma)\ge\sum_{i=1}^{n}\frac{1}{2}\min_{p^{i}\in P_{i}}\sum_{j=1}^{\left|p^{i}\right|}\wt{\left(p_{j}^{i}p_{j+1}^{i}\right)}\label{eq:LB}
\end{equation}
where, for each $i$, $P_{i}$ is the set of all walks of length $I_{i}(\pi,\sigma)$
starting from $\pi^{-1}(i)$ and ending at $\sigma^{-1}(i)$. Since
$\wtfn$ is decreasing, the minimum is attained by the walks $p^{i,\star}=(\pi^{-1}(i),\cdots,\ell_{i}-1,\ell_{i},\ell_{i}-1,\cdots,\sigma^{-1}(i))$
where $\ell_{i}$ is the solution to the equation 
\[
\ell_{i}-\pi^{-1}(i)+\ell_{i}-\sigma^{-1}(i)=I_{i}(\pi,\sigma)
\]
and thus $\ell_{i}=\left(\pi^{-1}(i)+\sigma^{-1}(i)+I_{i}(\pi,\sigma)\right)/2.$

We show next that there exists a transformation $\tau^{\star}$ such
that $p^{i,\tau^{\star}}=p^{i,\star}$ and thus equality in (\ref{eq:LB})
can be achieved. The transformation in question, $\tau^{\star}$,
transforms $\pi$ into $\sigma$ in $n$ rounds. In round $i$, $\tau^{\star}$
moves $i$ through a sequence of adjacent transpositions from position
$\pi^{-1}(i)$ to position $\sigma^{-1}(i)$. It can be seen that,
for each $i$, $p^{i,\tau}=(\pi^{-1}(i),\cdots,\ell_{i}'-1,\ell_{i}',\ell_{i}'-1,\cdots,\sigma^{-1}(i))$
for some $\ell_{i}'$. Since each transposition in $\tau$ decreases
the number of inversions by one, $\ell_{i}'$ also satisfies the equation
\[
\ell_{i}'-\pi^{-1}(i)+\ell_{i}'-\sigma^{-1}(i)=I_{i}(\pi,\sigma),
\]
implying that $\ell'_{i}=\ell_{i}$ and thus $p^{i,\tau^{\star}}=p^{i,\star}$.
Consequently, one has the following proposition.
\begin{prop}
For rankings $\pi,\sigma\in\mathbb{S}_{n}$, we have 
\[
\dist_{\wtfn}(\pi,\sigma)=\sum_{i=1}^{n}\frac{1}{2}\left(\sum_{j=\pi^{-1}(i)}^{\ell_{i}-1}\wt{(j\ j+1)}+\sum_{j=\sigma^{-1}(i)}^{\ell_{i}-1}\wt{(j\ j+1)}\right)
\]
where $\ell_{i}=\left(\pi^{-1}(i)+\sigma^{-1}(i)+I_{i}(\pi,\sigma)\right)/2$.\end{prop}
\begin{example}
Consider the rankings $\pi=4312$ and $e=1234$ and a decreasing weight
function $\wtfn$. We have $I_{i}(\pi,e)=2$ for $i=1,2$ and $I_{i}(\pi,e)=3$
for $i=3,4$. Furthermore, 
\begin{align*}
\ell_{1} & =\frac{3+1+2}{2}=3, & p^{1,\star} & =(3,2,1),\\
\ell_{2} & =\frac{4+2+2}{2}=4, & p^{2,\star} & =(4,3,2),\\
\ell_{3} & =\frac{2+3+3}{2}=4, & p^{3,\star} & =(2,3,4,3),\\
\ell_{4} & =\frac{1+4+3}{2}=4, & p^{4,\star} & =(1,2,3,4).
\end{align*}
The minimum weight transformation is 
\[
\tau^{\star}=\left(\underbrace{(32),(21)}_{1},\underbrace{(43),(32)}_{2},\underbrace{(43)}_{3}\right)
\]
 where the numbers under the braces are the element that is moved
by the indicated transpositions. The distance between $\pi$ and $e$
is 
\[
\dist_{\wtfn}(\pi,e)=\wt{(12)}+2\wt{(23)}+2\wt{(34)}.
\]

\end{example}

Note that the result above implies that at least for one class of
interesting weight functions that capture the importance of the position
in the ranking, the computation of the distance is of the same order
of complexity as that of standard Kendall's $\tau$ distance. Hence,
distance computation does not represent a bottleneck for the employment
of weighted distance metrics.

\section{Generalizing Kemeny's Approach}

We proceed by showing how Kemeny's axiomatic approach may be extended
further to introduce a number of new distances metrics useful in different
ranking scenarios.

The first distance applies when only certain subsets of transpositions
are allowed -- for example, when only elements of a class may be reordered
to obtain an aggregated ranking.
\begin{defn}
Consider a subset $G=\{g_{1},\cdots,g_{m}\}$ of $\mathbb{S}_{n}$
such that $g\in G$ implies that $g^{-1}\in G$. Rankings $\pi$ and
$\sigma$ are $G-$adjacent if there exist $g\in G$ such that $\pi=\sigma g$. 

A \emph{$G-$transformation} of $\pi$ into $\sigma$ is a vector
$(g_{1},\cdots,g_{k}),k\in\mathbb{N}$, with $g_{i}\in G,i\in[k]$,
such that $\sigma=\pi g_{1}g_{2}\cdots g_{k}$ where $k$ is the length
of the $G-$transformation. The set of $G-$transformations of $\pi$
into $\sigma$ is denoted by $A_{G}(\pi,\sigma)$. A \emph{minimum
$G-$transformation} is a $G-$transformation of minimum length.

Furthermore, $\omega$ is said to be $G-$between $\pi$ and $\sigma$
if there exists a minimal transformation $(g_{1},\cdots,g_{k})$ of
$\pi$ into $\sigma$ such that $\omega=\sigma g_{1}\cdots g_{j}$
for some $j\in[k]$.
\end{defn}

\begin{defn}
\label{def:uniform-distance}For a subset $G$ of $\mathbb{S}_{n}$,
a function $\dist:\mathbb{S}_{n}\to[0,\infty]$ is said to be a \emph{uniform
$G-$distance} if\end{defn}
\begin{enumerate}
\item $\dist$ is a metric.
\item $\dist$ is left-invariant.
\item For any $\pi,\sigma\in\mathbb{S}_{n}$, if $\omega$ is between $\pi$
and $\sigma$, then $\dist(\pi,\sigma)=\dist(\pi,\omega)+\dist(\omega,\sigma)$.
\item The smallest positive distance is one.
\end{enumerate}
Definition \ref{lem:on-the-line} also applies to $G-$betweenness
and can be restated as follows.
\begin{lem}
\label{lem:G-on-the-line} For a uniform $G-$distance $\dist$, and
for $\pi_{0},\cdots,\pi_{m}$ such that $\pi_{0}\dash\pi_{1}\dash\cdots\dash\pi_{m}$,
we have 
\[
\dist(\pi_{0},\pi_{m})=\sum_{k=1}^{m}\dist(\pi_{k-1},\pi_{k}).
\]
\end{lem}
\begin{rem}
\label{rem:uniform}For some choices of $G$, as in Lemma \ref{lem:two-ways}
and Lemma \ref{lem:T-two-ways} in the next section, one may show
that all elements of $G$ have distance one from the identity. For
such $G$, it is easy to see that the uniform $G-$distance $\dist$
exists and is unique, with 
\[
\dist(\pi,\sigma)=\min_{m}\left\{ m:(\tau_{1},\cdots,\tau_{m})\in A_{G}(\pi,\sigma)\right\} .
\]

\end{rem}

\begin{defn}
\label{def:weighted-distance}For a subset $G$ of $\mathbb{S}_{n}$,
a function $\dist:\mathbb{S}_{n}\to[0,\infty]$ is said to be a \emph{weighted
$G-$distance} if\end{defn}
\begin{enumerate}
\item $\dist$ is a pseudo-metric.
\item $\dist$ is left-invariant.
\item For any $\pi,\sigma\in\mathbb{S}_{n}$, if $\pi$ and $\sigma$ are
not $G-$adjacent, there exists a $\omega$ between $\pi$ and $\sigma$,
distinct from both, such that $\dist(\pi,\sigma)=\dist(\pi,\omega)+\dist(\omega,\sigma)$.\end{enumerate}
\begin{rem}
\label{rem:weighted}It is straightforward to see that the weighted
$G-$distance $\dist$ exists and is uniquely determined by the values
$\dist(g,e),g\in G$ as 
\[
\dist(\pi,\sigma)=\min_{\left(\tau_{1},\cdots,\tau_{m}\right)\in A_{G}(\pi,\sigma)}\sum_{i=1}^{m}\dist\left(\tau_{i},e\right)
\]
 where the minimum is taken over all $G-$transformations $\left(\tau_{1},\cdots,\tau_{m}\right)$
of $\pi$ into $\sigma$.
\end{rem}
As an example, let $G$ from Definitions \ref{def:uniform-distance}
and \ref{def:weighted-distance} be the set 
\[
\mathbb{T}_{n}=\{(ab):a,b\in[n],a\neq b\}
\]
of all transpositions.

The following lemma states that for a uniform $\mathbb{T}_{n}-$distance,
all transpositions have equal distance from identity.
\begin{lem}
\label{lem:T-two-ways}For a uniform $\mathbb{T}_{n}-$distance $\dist$,
we have 
\[
\dist\left((ab),e\right)=\dist\left((cd),e\right)
\]
for all transpositions $(ab)$ and $(cd)$.\end{lem}
\begin{IEEEproof}
For $\{a,b\}=\{c,d\}$, the lemma is obvious. We prove the lemma for
the case that $a,b,c,$ and $d$ are distinct. A similar argument
applies when $\{a,b\}$ and $\{c,d\}$ have one element in common.
The argument parallels that of Lemma \ref{lem:two-ways}. 

Let $\pi=(abcd)$, $\omega=(ad)\pi$, $\eta=(cd)\omega$ and note
that $e=(bc)\eta$. Since, $\pi\dash\omega\dash\eta\dash e$ by Lemma
\ref{lem:G-on-the-line} and left-invariance of $\dist$, we have
\begin{equation}
\dist(\pi,e)=\dist\left((ad),e\right)+\dist\left((cd),e\right)+\dist\left((bc),e\right).\label{eq:T-1}
\end{equation}
Similarly, let $\alpha=(bc)\pi$, $\beta=(ab)\alpha$, and note that
$e=(ad)\beta$. This shows
\begin{equation}
\dist(\pi,e)=\dist\left((bc),e\right)+\dist\left((ab),e\right)+\dist\left((ad),e\right).\label{eq:T-2}
\end{equation}
Equating the right-hand-sids of (\ref{eq:T-1}) and (\ref{eq:T-2})
yields $\dist\left((ab),e\right)=\dist\left((cd),e\right)$.
\end{IEEEproof}
By combining Remark \ref{rem:uniform} and Lemma \ref{lem:T-two-ways},
we arrive at the following theorem.
\begin{thm}
The uniform $\mathbb{T}_{n}-$distance exists and is unique. Namely,
\[
\dist(\pi,\sigma)=\min_{m}\left\{ m:(\tau_{1},\cdots,\tau_{m})\in A_{\mathbb{T}_{n}}(\pi,\sigma)\right\} ,
\]
(commonly known as Cayley's distance) is the unique $\mathbb{T}_{n}-$distance.
\end{thm}
Furthermore, Remark \ref{rem:weighted} leads to the following theorem.
\begin{thm}
\label{thm:T-weighted}The weighted $\mathbb{T}_{n}-$distance $\dist$
exists and is uniquely determined by the values $\dist(\tau,e)$,
$\tau\in\mathbb{T}_{n}$ as 
\[
\dist(\pi,\sigma)=\min_{\left(\tau_{1},\cdots,\tau_{m}\right)\in A_{\mathbb{T}_{n}}(\pi,\sigma)}\sum_{i=1}^{m}\dist\left(\tau_{i},e\right).
\]

\end{thm}
The weighted transposition distance can be used to model similarities
of objects in rankings wherein transposing two similar items induces
a smaller distance than transposing two dissimilar items \cite{farnoud2012sorting}.

\begin{rem}
Note that the generalization of Kemeny's axioms may also be applied
to arrive at a generalization of Borda's score-based rule. A step
in this direction was proposed by Young \cite{Young:1975fk}, who
showed that a set of axioms leads to a generalization of Borda's rule
wherein the $k$th preference of each ranking receives a score $s_{k}$,
not necessarily equal to $k$. This generalization of Borda's rule
may also be used to address the problem of top versus bottom in rankings.
In particular, one may assign Borda scores $s_{k}$ to the $k$th
preference with 
\begin{align*}
s_{k}=\sum_{l=1}^{k-1}\phi_{l} & ,
\end{align*}
 where $\phi_{k}$ is decreasing in $l$. For example, swapping two
elements at the top of the ranking of a given voter changes the scores
of each of the two corresponding objects by $\phi_{1}$ while a similar
swap at the bottom of the ranking, changes the scores by $\phi_{n-1}$.
Since $\phi_{1}\ge\phi_{n-1}$, changes to the top of the list, in
general, have a more significant affect on the aggregate ranking.
\end{rem}

\section{Distributed Vote Aggregation}

The novel distance metrics, scoring methods and underlying rank aggregation
problems discussed in the previous sections may be viewed as instances
rank aggregation of $m$ agents over a fully connected graph: i.e.\ every
agent has access to the ranking of all other agents and hence, fixing
the aggregation distance or scores and aggregation method (Kendall,
Borda,...) and assuming infinite computational power, each individual
can find an aggregate ranking of the society. Thus, assuming the uniqueness
of the aggregate ranking, agents come to a \textit{consensus} over
the aggregate ranking in one computational step. Nevertheless, one
can consider the more general problem of reaching consensus about
the aggregate ranking in an arbitrary network through local interactions.
In this section, we consider this problem over general networks and
provide an analysis of convergence for a specific choice of aggregation
method: i.e.\ the Borda aggregation method. The analysis of aggregation
methods for some other distance measures described in the paper is
postponed to the full version of the paper.

Let $G=([m],E)$ be a connected undirected graph over $m$ vertices
with the edge set $E$ that represents the connectivity pattern of
agents in a network%
\footnote{Many of the discussions in this section can be generalized for the
case of time-varying networks%
}. As before, we assume that each agent $i\in[m]$ has a ranking $\sigma_{i}$
over $n$ entities. There are multiple ways of distributed aggregation
of opinion in such a network, all of which are recursive schemes.

One way to perform distributed aggregation is through neighbor aggregation.
In this method, at discrete-time instances $t=0,1,\ldots$, each agent
maintains an estimate $\hat{\pi}_{i}(t)$ of the aggregate ranking.
At time $t$, each agent exchange its believe with his neighboring
agents. Then, at time $t+1$, agent $i$ sets its believe $\hat{\pi}_{i}(t+1)$
to be the aggregate ranking of all the estimates of the neighboring
agents at time $t$, including his own aggregation.

Another way to do distributed aggregation is through gossiping over
networks\cite{shah2009gossip}. Suppose that at each time instance
we pick an edge $\{i,j\}\in E$ with probability $p_{ij}>0$. Then,
agents $i$ and $j$ exchange their estimates $\hat{\pi}_{i}(t)$
and $\hat{\pi}_{j}(t)$ at time $t$ and they both let $\hat{\pi}_{i}(t+1)=\hat{\pi}_{j}(t+1)$
be the aggregation of $\hat{\pi}_{i}(t)$ and $\hat{\pi}_{j}(t)$.

\subsection{Gossiping Borda Vectors}

We describe next a distributed method using the Borda's scheme and
gossiping over networks. Let $b_{i}=b_{i}(0)$ be the vector of the
initial rankings of $n$ entities for agent $i$ (for Borda's method
we have the specific choice of $b_{i}=\pi_{i}^{-1}$). The goal is
to compute $\bar{b}=\frac{1}{m}\sum_{i=1}^{m}b_{i}(0)$. One immediate
solution to find $\bar{b}$ is through gossiping over the network
as described by the following algorithm:

\textit{Distributed Rank Aggregation:} 
\begin{enumerate}
\item At time $t\geq0$, pick an edge $\{i,i'\}\in E$ with probability
$P_{ii'}>0$ where $\sum_{\{i,i'\}\in E}P_{ii'}=1$, 
\item Let $i,i'$ exchange their estimate $b_{i}(t),b_{i'}(t)$ and let
$b_{i}(t+1)=b_{i'}(t+1)=\frac{1}{2}(b_{i}(t)+b_{i'}(t))$, 
\item For $\ell\not=i,i'$, let $b_{\ell}(t+1)=b_{\ell}(t)$. 
\end{enumerate}
As proven in \cite{Boyd06}, the above scheme approaches the average
as $t$ goes to infinity.
\begin{lem}
If $G=([m],E)$ is connected, then we almost surely have $\lim_{t\to\infty}b_{i}(t)=\bar{b}$. \end{lem}
\begin{IEEEproof}
The lemma is direct consequence of the results in \cite{Boyd06}.
\end{IEEEproof}
Note that in the distributed rank aggregation algorithm the ultimate
goal is to find the correct ordering of $\bar{b}=\frac{1}{m}\sum_{i=1}^{m}b_{i}(0)$
rather than the vector $\bar{b}$ itself. Thus, it is not important
that the estimates of the ranking vectors converges to $\bar{b}$,
but that the estimates of the actual ranks are correct. In other words,
if for some time $t$, for all agents $i$, the ordering of $b_{i}(t)$
matches the ordering of $\bar{b}$ for all $i\in[m]$, then the society
has already achieved consensus over the \textit{ranking} of the objects.
Here, we derive a probabilistic bound on the number of iterations
needed to probabilistically reach the optimum ranking.

Throughout the following discussions, without loss of generality we
may assume that $\bar{b}$ is ordered%
\footnote{Throughout this section we use superscript to denote the ranking of
objects.%
}, i.e.\ $\bar{b}^{1}\leq\bar{b}^{2}\leq\cdots\leq\bar{b}^{n}$. We
say that $t$ is a consensus time for the \textit{aggregate ranking}
if the ordering of $\bar{b}_{i}(t)$ matches the ordering of $\bar{b}$
for all $i\in[m]$. The following result follows immediately from
this definition: 
\begin{lem}
If $t$ is a consensus time for the ranking, then any $t'>t$ is a
consensus time for the ranking. \end{lem}
\begin{IEEEproof}
It suffice to show the result for $t'=t+1$. Let $\{i,i'\}$ be the
edge that is chosen randomly at time $t$. Since $t$ is a consensus
time for the ranking, we have $b_{i}^{1}(t)\leq\cdots\leq b_{i}^{n}(t)$
and $b_{i'}^{1}(t)\leq\cdots\leq b_{i'}^{n}(t)$, and thus we also
have 
\begin{align*}
b_{i}^{1}(t+1) & =\frac{1}{2}\left(b_{i}^{1}(t)+b_{i'}^{1}(t)\right)\\
 & \leq\cdots\leq b_{i}^{n}(t+1)\\
 & =\frac{1}{2}\left(b_{i}^{n}(t)+b_{i'}^{n}(t)\right),
\end{align*}
 which proves the claim. 
\end{IEEEproof}
Based on the lemma above, let us define \textit{the consensus time}
$T$ for the ordering to be: 
\[
T=\min\{t\geq0\mid\mbox{\ensuremath{t\:}is a consensus time for the ordering}.\}
\]

Note that for the random gossip scheme, $T$ is a random variable
and if we have an adapted process for the random choice of edges,
$T$ is a stopping time. Our goal is to provide a probabilistic bound
for $T$. For this, let $r^{j}=\min\left\{ \bar{b}^{j+1}-\bar{b}^{j},\bar{b}^{j}-\bar{b}^{j-1}\right\} $
and let $d^{j}=\max_{i}\bar{b}_{i}^{j}(0)-\min_{i}\bar{b}_{i}^{j}(0)$.
That is, $r_{j}$ is the minimum distance of the average rating of
$j$ from the neighboring objects and $d_{j}$ is the spread of the
initial ratings of the agents for the object $j$. Then, we have the
following result. 
\begin{thm}
For the consensus time $T$ of the ordering we have 
\[
P(T>t)\leq4m\lambda_{2}^{t}(W)\sum_{j=1}^{n}\left(\frac{d^{j}}{r^{j}}\right)^{2},
\]
 where $W=\sum_{\{i,i'\}\in E}P_{ii'}\left(I-\frac{1}{2}(e_{i}-e_{i'})(e_{i}-e_{i'})^{T}\right)$,
$e_{i}=[0\ \cdots\ 0\ 1\ 0\ \cdots\ 0]^{T}$ is an $m\times1$ vector
with $i$th element equal to one, and $\lambda_{2}(W)$ is the second
largest eigenvalue of $W$.\end{thm}
\begin{IEEEproof}
Let $b^{j}(t)$ be the vector obtained by the rating of the $m$ agents
at time $t$ for object $j$ and let $y^{j}(t)=b^{j}(t)-\bar{b}^{j}$.
Note that if $\|y^{j}(t)\|^{2}\leq\left(\frac{r^{j}}{2}\right)^{2}$,
then this means that $|b_{i}^{j}(t)-\bar{b}^{j}|\leq\frac{r^{j}}{2}$
for all $i$. Thus, if for all $j\in[n]$ we have $\|y^{j}(t)\|^{2}\leq\left(\frac{r^{j}}{2}\right)^{2}$,
then it follows that: 
\[
b_{i}^{j}(t)\leq\bar{b}^{j}+\frac{r^{j}}{2}\leq\frac{1}{2}(\bar{b}^{j}+\bar{b}^{j+1}),
\]
 where the last inequality follows from the fact that $r^{j}\leq\bar{b}^{j+1}-\bar{b}^{j}$.
Similarly, we have: 
\[
b_{i}^{j+1}(t)\geq\bar{b}^{j+1}-\frac{r^{j+1}}{2}\geq\frac{1}{2}(\bar{b}^{j+1}+\bar{b}^{j}),
\]
 which follows from $r^{j+1}\leq\bar{b}^{j+1}-\bar{b}^{j}$. Hence,
we have 
\begin{align*}
b_{i}^{1}(t) & \leq\frac{1}{2}(\bar{b}^{2}+\bar{b}^{1})\leq b_{i}^{2}(t)\leq\frac{1}{2}(\bar{b}^{3}+\bar{b}^{2})\\
 & \leq\cdots\leq\frac{1}{2}(\bar{b}^{n-1}+\bar{b}^{n})\leq b_{i}^{m}(t),
\end{align*}
 and so $t$ is a consensus time for the algorithm. Thus, 
\[
\{T>t\}\subseteq\bigcup_{j=1}^{n}\left\{ \|y^{j}(t)\|^{2}\geq\left(\frac{r^{j}}{2}\right)^{2}\right\} 
\]
 and hence, using the union bound, we obtain 
\begin{equation}
P(T>t)\leq\sum_{j=1}^{n}P\left(\|y^{j}(t)\|^{2}\geq\left(\frac{r^{j}}{2}\right)^{2}\right).\label{eqn:unionbound}
\end{equation}
 Markov's inequality implies that 
\[
P\left(\|y^{j}(t)\|^{2}\geq\left(\frac{r^{j}}{2}\right)^{2}\right)\leq\left(\frac{2}{r^{j}}\right)^{2}E\left[\|y^{j}(t)\|^{2}\right].
\]

Using the analysis in \cite{Boyd06}, it can be shown that 
\[
E\left[\|y^{j}(t)\|^{2}\right]\leq\lambda_{2}^{t}(W)\|y^{j}(0)\|^{2}\leq m\left(d^{j}\right)^{2}.
\]
 Combining the above two relations, we find 
\[
P\left(\|y^{j}(t)\|^{2}\geq\left(\frac{r^{j}}{2}\right)^{2}\right)\le4m\lambda_{2}^{t}\left(\frac{d^{j}}{r^{j}}\right)^{2}.
\]
 Replacing the last inequality in (\ref{eqn:unionbound}), proves
the assertion.
\end{IEEEproof}
Note that from \cite{Boyd06}, if $G$ is connected, then we have
$\lambda_{2}<1$ and thus the probability $P(T>t)$ decays exponentially.
\bibliographystyle{ieeetr}
\bibliography{IEEEfull,bib}

\begin{thebibliography}{10}

\bibitem{kemeney1959mathematics}
J.~G. Kemeny, ``Mathematics without numbers,'' {\em Daedalus}, vol.~88, no.~4,
  pp.~pp. 577--591, 1959.

\bibitem{cook1985ordinal}
W.~D. Cook and M.~Kress, ``Ordinal ranking with intensity of preference,'' {\em
  Management Science}, vol.~31, pp.~26--32, 01 1985.

\bibitem{dwork2001rank}
C.~Dwork, R.~Kumar, M.~Naor, and D.~Sivakumar, ``Rank aggregation revisited,''
  {\em Manuscript (Available at:
  www.eecs.harvard.edu/~michaelm/CS222/rank2.pdf)}, 2001.

\bibitem{dwork2001rank-web}
C.~Dwork, R.~Kumar, M.~Naor, and D.~Sivakumar, ``Rank aggregation methods for
  the web,'' in {\em Proceedings of the 10th international conference on World
  Wide Web}, pp.~613--622, ACM, 2001.

\bibitem{sculley2007rank}
D.~Sculley, ``Rank aggregation for similar items,'' in {\em Proceedings of the
  Seventh SIAM International Conference on Data Mining}, Citeseer, 2007.

\bibitem{schalekamp2009rank}
F.~Schalekamp and A.~van Zuylen, ``Rank aggregation: Together we're strong,''
  {\em Proc. of 11th ALENEX}, pp.~38--51, 2009.

\bibitem{kumar2010gdr}
R.~Kumar and S.~Vassilvitskii, ``Generalized distances between rankings,'' in
  {\em Proceedings of the 19th international conference on World wide web}, WWW
  '10, (New York, NY, USA), pp.~571--580, ACM, 2010.

\bibitem{diaconis1988group}
P.~Diaconis, ``Group representations in probability and statistics,'' {\em
  Lecture Notes-Monograph Series}, vol.~11, 1988.

\bibitem{arrow1963social}
K.~J. Arrow, {\em Social choice and individual values}.
\newblock No.~12, Yale Univ Pr, 1963.

\bibitem{borda1784}
J.-C. de~Borda, ``M\'{e}moire sur les \'{e}lections au scrutin,'' {\em Histoire
  de l'Acad\'{e}mie royale des sciences}, 1784.

\bibitem{young1974-an-axiomatization}
H.~P. Young, ``An axiomatization of borda's rule,'' {\em Journal of Economic
  Theory}, vol.~9, no.~1, pp.~43 -- 52, 1974.

\bibitem{farnoud2012sorting}
F.~Farnoud and O.~Milenkovic, ``Sorting of permutations by cost-constrained
  transpositions,'' {\em Information Theory, IEEE Transactions on}, vol.~58,
  pp.~3 --23, Jan. 2012.

\bibitem{Young:1975fk}
H.~Young, ``Social choice scoring functions,'' {\em SIAM Journal on Applied
  Mathematics}, vol.~28, pp.~824--838, 06 1975.

\bibitem{shah2009gossip}
D.~Shah, {\em Gossip Algorithms}, vol.~3.
\newblock Foundations and Trends in Networking, 2009.

\bibitem{Boyd06}
S.~Boyd, A.~Ghosh, B.~Prabhakar, and D.~Shah, ``Randomized gossip algorithms,''
  {\em IEEE Transactions on Information Theory}, vol.~52, no.~6,
  pp.~2508--2530, 2006.

\end{thebibliography}

\end{document}